\documentclass[12pt]{article}
\usepackage{amssymb}

\usepackage{amsmath}
\usepackage{graphicx}

\newtheorem{corollary}{Corollary}
\newtheorem{definition}{Definition}
\newtheorem{example}{Example}
\newtheorem{lemma}{Lemma}

\newtheorem{theorem}{Theorem}
\newenvironment{proof}[1][Proof]{\textbf{#1.} }{\ \rule{0.5em}{0.5em}}

\begin{document}

\title{Dedekind order completion of $C(X)$ by Hausdorff continuous functions}
%running title:Dedekind order completion of $C(X)$
\author{R Anguelov \\
%EndAName
Department of Mathematics and Applied Mathematics\\
University of Pretoria\\
Pretoria 0002\\
anguelov@scientia.up.ac.za}
\date{}
\maketitle

\begin{abstract}
The concept of Hausdorff continuous interval valued functions, developed
within the theory of Hausdorff approximations and originaly defined for
interval valued functions of one real variable is extended to interval
valued functions defined on a topological space $X$. The main result is that
the set $\mathbb{H}_{f\!t}(X)$ of all finite Hausdorff continuous functions
on any topological space $X$ is Dedekind order complete. Hence it contains
the Dedekind order completion of the set $C(X)$ of all continuous real
functions defined on $X$ as well as the Dedekind order completion of the set
$C_{b}(X)$ of all bounded continuous functions on $X$. Under some general
assumptions about the topological space $X$ the Dedekind order completions
of both $C(X)$ and $C_{b}(X)$ are characterised as subsets of $\mathbb{H}%
_{f\!t}(X)$. This solves a long outstanding open problem about the Dedekind
order completion of $C(X)$. In addition, it has major applications to the
regularity of solutions of large classes of nonlinear PDEs.
\end{abstract}

2000 Mathematics Subject Classification 35F20, 54C30, 26E25

\section{Introduction}

The fact that the set $C(X)$ of all continuous real valued function on a
topological space $X$ is generally not Dedekind order complete with respect
to the point-wise defined partial order%
\begin{equation}
f\leq g\Longleftrightarrow f(x)\leq g(x),\;x\in X.  \label{order}
\end{equation}%
is well known and can be shown by trivial examples. We consider the problem
of constructing a Dedekind order completion of $C(X)$ through functions
defined on the same space $X$. An earlier result by Dilworth, see \cite%
{Dilworth}, gives a Dedekind order completion of the set $C_{b}(X)$ of all
bounded continuous functions on a completely regular topological space\ $X$
through the so called normal upper semi-continuous functions on $X$, the
Dedekind order completion of $C(X)$ in the general case when $C(X)$ contains
unbounded functions remaining an open problem. Here we obtain Dedekind order
completions of both $C_{b}(X)$ and $C(X)$ through Hausdorff continuous
interval valued functions. The significance of this result is emphasized by
the fact that important applications involve sets of continuous functions
which are not bounded. In \cite{Rosinger} it was shown that arbitrary
nonlinear PDEs defined by continuous, not necessarily smooth or analytic
expressions, have solutions that can be assimilated with Lebesgue measurable
functions. This powerful earlier existence results can now significantly be
improved with respect to the regularity of solutions by showing that the
solutions are in fact Hausdorff continuous, see for details \cite%
{Anguelov-Rosinger}.

We recall that a partially ordered set $P$ is called Dedekind order complete
if every subset of $P$ which is bounded from above has a supremum in $P$ and
every subset of $P$ which is bounded from below has an infimum in $P$. A
general result on the Dedekind order completion of partially ordered sets
was established by MacNeilly in 1937 (see \cite{Luxemburg} for a more recent
presentation). The problem of the order completion of $C(X)$ is particularly
addressed in \cite{Mack}. As stated above, the problem of constructing a
Dedekind order completion of $C(X)$ as a set of functions on $X$ was
partially addressed by Dilworth.\ More precisely, it was proved in \cite%
{Dilworth} that if $X$ is completely regular, then the Dedekind order
completion of $C_{b}(X)$ is isomorphic with the lattice of the normal upper
semi-continuous functions defined on $X$.\ Our approach is to consider $C(X)$
as a subset of the set of interval valued functions defined on $X$ and find
the Dedekind order completion of $C(X)$ within this set.

Denote by $\mathbb{I\,}\overline{\mathbb{R}}$ the set of all usual or
extended real intervals
\begin{equation*}
\mathbb{I\,}\overline{\mathbb{R}}=\{[\underline{a},\overline{a}]:\underline{a%
},\overline{a}\in \overline{\mathbb{R}}=\mathbb{R\cup \{\pm \infty \}},\;%
\underline{a}\leq \overline{a}\}.
\end{equation*}%
We consider on $\mathbb{I\,}\overline{\mathbb{R}}$ the partial order $\leq $
defined in \cite{Markov} through
\begin{equation}
\lbrack \underline{a},\overline{a}]\leq \lbrack \underline{b},\overline{b}%
]\Longleftrightarrow \underline{a}\leq \underline{b},\;\overline{a}\leq
\overline{b}.  \label{iorder}
\end{equation}

In the set of interval valued functions $\mathbb{A}(X)=\{f:X\rightarrow
\mathbb{I\,}\overline{\mathbb{R}}\}$, a partial order is induced by (\ref%
{iorder}) in a point-wise way similar with (\ref{order}). Identifying $x\in
\overline{\mathbb{R}}$ with $[x,x]\in \mathbb{I\,}\overline{\mathbb{R}}$ we
consider $\overline{\mathbb{R}}$ as a subset of $\mathbb{I\,}\overline{%
\mathbb{R}}$. In this way $\mathbb{A}(X)$ contains the set of extended real
valued functions $\mathcal{A}(X)=\{f:X\rightarrow \overline{\mathbb{R}}\}.$
Hence%
\begin{equation*}
C(X)\subset \mathcal{A}(X)\subset \mathbb{A}(X).
\end{equation*}%
In Section 2 we define the set $\mathbb{H}_{f\!t}(X)$ of all the Hausdorff
continuous interval valued functions with values finite real intervals, in
which case we shall have%
\begin{equation}
C(X)\subset \mathbb{H}_{f\!t}(X)\subset \mathbb{A}(X).  \label{incl1}
\end{equation}%
In Sections 3-5 the following results are presented:

\begin{enumerate}
\item The set $\mathbb{H}_{f\!t}(X)$ is Dedekind order complete. Hence it
contains the Dedekind order completions $C(X)^{\#}$ and $C_{b}(X)^{\#}$ of $%
C(X)$ and $C_{b}(X)$, respectively.

\item For a completely regular topological space $X$
\begin{equation*}
C_{b}(X)^{\#}=\mathbb{H}_{b}(X)\subseteq\mathbb{H}_{f\!t}(X)
\end{equation*}
where $\mathbb{H}_{b}(X)$ is the set of all bounded Hausdorff continuous
functions defined on $X$.

\item For a completely regular topological space $X$%
\begin{equation*}
C(X)^{\#}=\mathbb{H}_{cm}(X)\subseteq\mathbb{H}_{f\!t}(X)
\end{equation*}
where $\mathbb{H}_{cm}(X)$ is the set of all Hausdorff continuous functions
with continuous majorant and continuous minorant on $X$.

\item If $X$ is a metric space then%
\begin{equation*}
C(X)^{\#}=\mathbb{H}_{cm}(X)=\mathbb{H}_{f\!t}(X)
\end{equation*}
\end{enumerate}

The characterization of the topological spaces $X$ for which $\mathbb{H}%
_{cm}(X)=\mathbb{H}_{f\!t}(X)$ remains an open problem.

Historically, the interval analysis, or the analysis of interval valued
functions, is associated with, so called, validated computing where
algorithms generating validated bounds for the exact solutions of
mathematical problems are designed and investigated \cite{Alefeld}, \cite%
{Kramer}. However, interest in the interval valued functions comes also from
other branches of mathematics such as nonlinear partial differential
equations, see \cite{Anguelov-Rosinger} which strengthens the results in %
\cite{Rosinger}, and approximation theory \cite{Sendov}. In fact, Hausdorff
continuous functions of one real variable were first introduced by Sendov %
\cite{Sendov} in connection with Hausdorff approximations of real functions
of real argument. The name Hausdorff continuous is due to the
characterization of these functions in terms of the Hausdorff distance
between the graphs of real functions as defined in \cite{Sendov}. The
concept was further developed in \cite{Anguelov-Markov} as part of the
analysis of interval valued functions. Here, as a new departure in applying
the ideas of interval analysis, we consider the concept of Hausdorff
continuity for interval valued functions defined on a topological space $X$
and apply it to the Dedekind order completion of $C_{b}(X)$ and $C(X)$.

\section{Hausdorff Continuous Interval Valued Functions}

For every $x\in X,$ $\mathcal{V}_{x}$ denotes the set of neighborhoods of $x$%
. We consider, \cite{Baire}, the pair of mappings $I:\mathbb{A}%
(X)\rightarrow \mathcal{A}(X)$, $S:\mathbb{A}(X)\rightarrow \mathcal{A}(X)$,
called lower Baire, and upper Baire operators, respectively, where for every
function $f\in \mathbb{A}(X)$ and $x\in X$, we have%
\begin{align}
I(f)(x)& =\sup_{V\in \mathcal{V}_{x}}\inf \{z\in f(y):y\in V\},  \label{lbf}
\\
S(f)(x)& =\inf_{V\in \mathcal{V}_{x}}\sup \{z\in f(y):y\in V\}.  \label{ubf}
\end{align}%
The operator $F:\mathbb{A}(X)\rightarrow \mathbb{A}(X)$ defined by
\begin{equation*}
F(f)(x)=[I(f)(x),S(f)(x)],\text{ }f\in \mathbb{A}(X),\text{ }x\in X,
\end{equation*}%
is called graph completion.

Let us note that the lower Baire operator $I:f\rightarrow I(f)$, the upper
Baire operator $S:f\rightarrow S(f)$ and the graph completion operator $%
F:f\rightarrow F(f)=[I(f),S(f)]$ are all monotone with respect to the order $%
\leq$ in $\mathbb{A}(X)$, which means that for every two functions $f,g\in%
\mathbb{A}(X)$ we have%
\begin{equation}
f\leq g\Longrightarrow I(f)\leq I(g)\text{, }S(f)\leq S(g)\text{, }F(f)\leq
F(g).  \label{opres}
\end{equation}
The operator $F$ is also monotone about inclusion
\begin{equation*}
f(x)\subseteq g(x),x\in X\Longrightarrow F(f)(x)\subseteq F(g)(x)\text{, }%
x\in X.
\end{equation*}
Furthermore, all three operators are idempotent, i.e. for every $f\in
\mathbb{A}(X)$%
\begin{equation}
I(I(f))=I(f),\text{ }S(S(f))=S(f),\text{ }F(F(f))=F(f).  \label{idemp}
\end{equation}

The fixed points of the operators $I$ and $S$ are the lower and upper
semi-continuous functions, respectively, defined on $X$. Let us recall the
definitions of lower and upper semi-continuity.\cite{Baire}

\begin{definition}
A function $f:X\rightarrow\overline{\mathbb{R}}$ is called lower
semi-continuous at $x\in X$ if for every $m<f(x)$ there exists $V\in
\mathcal{V}_{x}$ such that $m<f(y)$ for all $y\in V$. If $f(x)=-\infty$,
then $f$ is assumed lower semi-continuous at $x$.
\end{definition}

\begin{definition}
A function $f:X\rightarrow\overline{\mathbb{R}}$ is called upper
semi-continuous at $x\in X$ if for every $m>f(x)$ there exists $V\in
\mathcal{V}_{x}$ such that $m>f(y)$ for all $y\in V$. If $f(x)=+\infty$,
then $f$ is assumed upper semi-continuous at $x$.
\end{definition}

\begin{definition}
A function $f:X\rightarrow\overline{\mathbb{R}}$ is called lower (upper)
semi-continuous on $X$ if it is lower (upper) semi-continuous at every point
of $X$.
\end{definition}

It is easy to see that for every $f\in\mathbb{A}(X)$ the functions $I(f)$
and $S(f)$ are, respectively, lower and upper semi-continuous on $X.$
Furthermore, if $f\in\mathcal{A}(X)$ then%
\begin{align*}
f\text{ - lower semi-continuous on }X & \Longleftrightarrow I(f)=f, \\
f\text{ - upper semi-continuous on }X & \Longleftrightarrow S(f)=f.
\end{align*}

\begin{definition}
\label{dhcon}A function $f\in\mathbb{A}(X)$ is called Hausdorff continuous,
or H-continuous, if for every function $g\in\mathbb{A}(X)$ which satisfies
the inclusion $g(x)\subseteq f(x),$ $x\in X$, we have $F(g)(x)=f(x),$ $x\in
X $.
\end{definition}

We denote by $\mathbb{H}(X)$ the subset of $\mathbb{A}(X)$ consisting of all
the H-continuous functions while by $\mathbb{H}_{f\!t}(X)$ we denote the set
of all H-continuous functions on $X$ which assume finite values for every $%
x\in X$, that is,%
\begin{equation*}
\mathbb{H}_{f\!t}(X)=\{f\in\mathbb{H}(X):f(x)\subseteq\mathbb{R},\text{ }%
x\in X\}
\end{equation*}
It is easy to see that all usual continuous real-valued functions defined on
$X$ are H-continuous. Indeed, since a continuous function $f$ is both lower
and upper semi-continuous we have
\begin{equation*}
F(f)=[I(f),S(f)]=[f,f]=f.
\end{equation*}
Furthermore, if $g\in\mathbb{A}(X)$ is such that $g(x)\subseteq f(x),$ $x\in
X,$ then $g(x)=f(x),$ $x\in X,$ because $f(x)$ is a point interval for all $%
x\in X$. Hence $F(g)(x)=F(f)(x)=f(x),$ $x\in X$ and the H-continuity of $f$
follows from the Definition \ref{dhcon} above. Thus $C(X)\subseteq \mathbb{H}%
(X).$ Moreover, since the functions in $C(X)$ assume values which are finite
real numbers we have%
\begin{equation*}
C(X)\subseteq\mathbb{H}_{f\!t}(X)\subseteq\mathbb{H}(X).
\end{equation*}
We will show in the next section that the set $\mathbb{H}_{f\!t}(X)$ is
Dedekind order complete. Therefore, in view of the above inclusion, for the
characterization of the Dedekind order completion of $C(X)$ we only need to
consider $\mathbb{H}_{f\!t}(X)$.

Let us note that the set $\mathbb{H}_{f\!t}(X)$ is certainly wider than $%
C(X) $. Here are some examples of H-continuous functions which are \emph{not}
continuous.

\begin{example}
Let $X=\mathbb{R}$. For $x\in X$%
\begin{equation*}
s(x)=\left\{
\begin{array}{ccc}
1 & \text{if} & x>0 \\
\lbrack-1,1] & \text{if} & x=0 \\
-1 & \text{if} & x<0%
\end{array}
\right. .
\end{equation*}
\end{example}

\begin{example}
Let $X=\mathbb{R}^{2}$. For $x\in X$%
\begin{equation*}
f(x)=\left\{
\begin{array}{ccc}
s\left( \sin\frac{1}{||x||}\right) & \text{if} & x\neq(0,0) \\
\lbrack-1,1] & \text{if} & x=(0,0)%
\end{array}
\right. .
\end{equation*}
\end{example}

\begin{example}
\label{shock} Consider a typical example of solution of a nonlinear shock
wave equation. Here $X=\{(x,t):t\geq 0\}\subseteq \mathbb{R}^{2}$ and for $%
x\in X$\vspace{3mm} \newline
\begin{minipage}{6.5cm}
$y(x,t)\!\!=\!\!\left\{ \!\!%
\begin{array}{cll}
1\bigskip & \!\text{if} & \!\!t\!\in \![0,1),x<t\!-\!1 \\
\frac{x}{t-1}\bigskip & \!\text{if} & \!\!t\!\in \![0,1),x\!\in
\![t\!-\!1,\!0] \\
0\bigskip & \!\text{if} & \!\!t\!\in \![0,1),x>0 \\
1\bigskip & \!\text{if} & \!\!t\!\geq 1\text{ , \ }x<\frac{t-1}{2} \\
\lbrack \!-\!1,\!1]\bigskip & \!\text{if} & \!\!t\!\geq 1\text{ , \ }x=\frac{%
t-1}{2} \\
0 & \!\text{if} & \!\!t\!\geq 1\text{ , \ }x>\frac{t-1}{2}%
\end{array}%
\right.$\end{minipage} \hspace{-0.5cm}
\begin{minipage}{8cm}
\includegraphics[height=160pt, width=170pt]{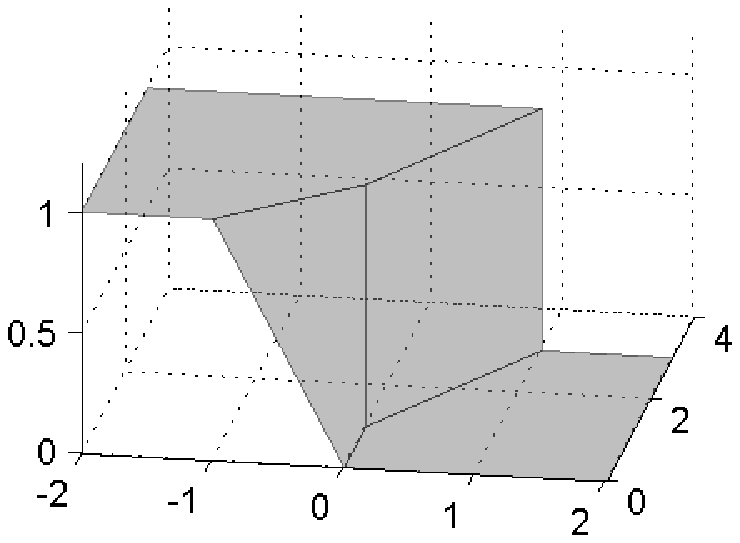}%
\end{minipage}

\end{example}

Let $f\in \mathbb{A}(X).$ For every $x\in X$ the value of $f$ is an interval
$[\underline{f}(x),\overline{f}(x)]$. Hence, the function $f$ can be written
in the form $f=[\underline{f},\overline{f}]$ where $\underline{f},\overline{f%
}\in \mathcal{A}(X)$ and $\underline{f}\leq \overline{f}.$ The lower and
upper Baire operators of the interval valued function $f$ can be
conveniently represented in terms of the functions $\underline{f}$ and $%
\overline{f}$. Indeed, from (\ref{lbf}) and (\ref{ubf}) it is easy to see
that
\begin{equation}
I(f)=I(\underline{f})\text{ \ and \ }S(f)=S(\overline{f}).  \label{ifsf}
\end{equation}%
Hence $F(f)$ can be written in the form%
\begin{equation}
F(f)=[I(\underline{f}),S(\overline{f})].  \label{fff}
\end{equation}%
Therefore,
\begin{equation}
F(f)=f\Longleftrightarrow \underline{f}=I(\underline{f}),\text{ }\overline{f}%
=S(\overline{f})\Longleftrightarrow \left\{
\begin{tabular}{l}
$\overline{f}$ - upper semi-continuous \\
$\underline{f}\text{ - lower semi-continuous}$%
\end{tabular}%
\ \ \right. .  \label{fscont}
\end{equation}

It is easy to see that if $f$ is H-continuous then it is a fixed point of
the operator $F$. Indeed, from the Definition \ref{dhcon} and the inclusion $%
f(x)\subseteq f(x),$ $x\in X,$ it follows that $F(f)=f$. Thus, in view of (%
\ref{fscont}), we have that if $f=[\underline{f},\overline{f}]$ is
H-continuous then the functions $\underline{f},$ $\overline{f}$ are lower
and upper semi-continuous functions, respectively. The following theorem
provides convenient criteria for recognizing H-continuous functions. These
criteria are discussed in \cite{Sendov}. However, since there only the case $%
X\subseteq \mathbb{R}$ is considered, we provide here a short proof.

\begin{theorem}
\label{Hcondi}Let $f=[\underline{f},\overline{f}]\in\mathbb{A}(X)$. The
following conditions are equivalent\newline
\hspace*{1cm}a) the function $f$ is H-continuous\newline
\hspace*{1cm}b) $F(\underline{f})=F(\overline{f})=f$\newline
\hspace*{1cm}c) $S(\underline{f})=\overline{f},$ $I(\overline {f})=%
\underline{f}$
\end{theorem}

\begin{proof}
a)$\Longrightarrow$b): Since $f$ is H-continuous using the Definition \ref%
{dhcon} we obtain
\begin{equation*}
\underline{f}(x)\subseteq f(x),x\in X\Longrightarrow F(\underline{f}%
)(x)=F(f)(x)=f(x),x\in X
\end{equation*}
In a similar way we prove $F(\overline{f})=f.$

b)$\Longrightarrow$c): We have%
\begin{align*}
\lbrack I(\underline{f}),S(\underline{f})] & =F(\underline{f}))= f=[%
\underline{f},\overline{f}]\Longrightarrow S(\underline{f} )=\overline {f},
\\
\lbrack I(\overline{f}),S(\overline{f})] & =F(\overline{f}))=f=[\underline {f%
},\overline{f} ]\Longrightarrow I(\overline{f})=\underline{f}
\end{align*}

c)$\Longrightarrow$a):\ Let us assume that $g=[\underline{g},\overline{g}
]\in\mathbb{A}(X)$ is such that $g(x)\subseteq f(x),x\in X$. Using the
monotonicity of the operators $I$ and $S$, thus from the inequalities
\begin{equation*}
\underline{f} \leq\underline{g}\leq\overline{g}\leq\overline{f}
\end{equation*}
it follows that%
\begin{align*}
S(\underline{f} ) & \leq S(\overline{g})\leq S(\overline{f})=\overline{f} ,
\\
\underline{f} & =I(\underline{f})\leq I(\underline{g} )\leq I(\overline{f}).
\end{align*}

The equalities in c) imply that $S(\overline{g})=\overline{f} $ and $I(%
\underline{g})=\underline{f}$. Therefore, $F(g)=f$, which shows that $f$ is
H-continuous.
\end{proof}

\begin{theorem}
\label{FSI}Let $f\in\mathbb{A}(X)$. Both functions $F(S(I(f)))$ and $%
F(I(S(f)))$ are H-continuous.
\end{theorem}

\begin{proof}
Denote $g=[\underline{g},\overline{g} ]=F(S(I(f))$. Clearly the functions $%
\underline{g}=I(S(I(f)))$ and $\overline{g}=S(I(f))]$ are lower and upper
semi-continuous, respectively. Using the monotonicity of the operators $I$
and $S$ we obtain the inequalities%
\begin{align*}
S(\underline{g} ) & =S(I(S(I(f))))\geq S(I(I(f)))=S(I(f))=\overline{g} , \\
S(\underline{g}) & \leq S(\overline{g})=\overline{g},
\end{align*}
which imply $S(\underline{g})=\overline{g}$. Furthermore, we have%
\begin{equation*}
I(\overline{g})=I(S(I(f)))=\underline{g}.
\end{equation*}
Then, it follows from Theorem \ref{Hcondi} (a and c) that function $%
g=F(S(I(f))$ is H-continuous. The H-continuity of $F(I(S(f)))$ is proved in
a similar way.
\end{proof}

\begin{theorem}
\label{tcont}Let $f=[\underline{f},\overline{f}]$ be an H-continuous
function on $X$. \newline
a) If $\underline{f}$ or $\overline{f}$ is continuous at a point $a\in X$
then $\underline{f}(a)=\overline{f}(a)$.\newline
b) If $\underline{f}(a)=\overline{f}(a)$ for some $a\in X$ then both $%
\underline{f}$ and $\overline{f}$ are continuous at $a$.
\end{theorem}

\begin{proof}
a) Let $\underline{f}$ be continuous at $a\in X$. Assume that the equality
at a) does not hold. Thus, $\underline{f}(a)<\overline{f}(a)$. Then, since $%
\underline{f}$ is lower semi-continuous at $a$ there exists $V_{1} \in%
\mathcal{V}_{a}$ such that $\underline{f}(y)<\frac{1}{2}(\underline{f}(a)+%
\overline{f} (a))$ whenever $y\in V_{1}$. Hence
\begin{align*}
S(\underline{f})(a) & =\inf_{V\in\mathcal{V}_{a}}\sup\{\underline{f}
(y):y\in V\}\leq\sup\{\underline{f}(y):y\in V_{1} \} \\
& \leq\frac{1}{2}(\underline{f}(a)+\overline{f}(a))<\overline{f}(a).
\end{align*}
This, according to Theorem \ref{Hcondi} (a and c), implies that $f$ is not
H-continuous which contradicts the condition of the theorem. Thus, $%
\underline{f}(a)=\overline{f}(a)$. The case when $\overline{f} $ is
continuous is treated in the same way.

b) Let $\varepsilon>0$. Using the semi-continuity of $\underline{f}$ and $%
\overline{f}$ we obtain that there exists a $V_{2}\in\mathcal{V}_{a}$ such
that for every $y\in V_{2}$ we have
\begin{align*}
\overline{f}(y) & \geq\underline{f} (y)>\underline{f}(a)-\varepsilon
=f(a)-\varepsilon, \\
\underline{f}(y) & \leq\overline{f} (y)<\overline{f}(a)+\varepsilon
=f(a)+\varepsilon.
\end{align*}
Therefore, $|\underline{f}(y)-\underline{f}(a)|=|\underline{f}%
(y)-f(a)|<\varepsilon$ and $|\overline{f}(y)-\overline{f}(a)|=|\overline {f}%
(y)-f(a)|<\varepsilon$ for all $y\in V_{2}$, which implies that both $%
\underline{f}$ and $\overline{f}$ are continuous at $a$.
\end{proof}

The following theorem shows a similarity between continuous functions and
H-continuous functions.

\begin{theorem}
\label{tindent}Let $f,g$ be H-continuous on $X$ and let $D$ be a dense
subset of $X$. Then%
\begin{align*}
\text{a) \ \ }f(x) & \leq g(x),\text{ }x\in D\Longrightarrow f(x)\leq g(x),%
\text{ }x\in X,\text{ \ \ \ \ \ \ \ \ \ \ \ \ \ \ \ \ \ \ \ \ } \\
\text{b) \ \ }f(x) & =g(x),\text{ }x\in D\Longrightarrow f(x)=g(x),\text{ }%
x\in X.\text{ \ \ \ \ \ \ \ \ \ \ \ \ \ \ \ \ \ \ \ \ }
\end{align*}
\end{theorem}

\begin{proof}
a) Let $f=[\underline{f},\overline{f}]$ and let $g=[\underline{g},\overline {%
g}]$. Let also $x\in X$ and $V_{1},V_{2}$ be arbitrary neighborhoods of $x$.
Since $D$ is dense in $X$ there exists $z_{0}\in D\cap(V_{1}\cap V_{2})$,
that is $z_{0}\in D\cap V_{1}$ and $z_{0}\in D\cap V_{2}$. Therefore%
\begin{equation*}
\inf\{\underline{f}(y):y\in V_{1}\}\leq\underline{f}(z_{0})\leq\underline {g}%
(z_{0})\leq\overline{g}(z_{0})\leq\sup\{\overline{g}(y):y\in V_{2}\}
\end{equation*}
Using that $V_{1}$ and $V_{2}$ are chosen independently we have%
\begin{equation*}
I(\underline{f})(x)=\sup_{V_{1}\in\mathcal{V}_{x}}\inf\{\underline{f}%
(y):y\in V_{1}\}\leq\inf_{V_{2}\in\mathcal{V}_{x}}\sup\{\overline{g}(y):y\in
V_{2}\}=S(\overline{g})(x).
\end{equation*}
Since\ $\underline{f}$ and $\overline{g}$ are lower and upper
semi-continuous functions, respectively, the above inequality implies%
\begin{equation}
\underline{f}(x)\leq\overline{g}(x),\text{ }x\in X.  \label{ineq1}
\end{equation}

Furthermore, $S(\underline{f})=\overline{f}$ and $I(\overline{g})=\underline{%
g}$, because both $f$ and $g$ are H-continuous, see Theorem \ref{Hcondi}.
Hence, using the monotonicity of the operators $I$ and $S$, see (\ref{opres}%
), and inequality (\ref{ineq1}) we obtain%
\begin{align*}
\overline{f}& =S(\underline{f})\leq S(\overline{g})=\overline{g}, \\
\underline{f}& =I(\underline{f})\leq I(\overline{g})=\underline{g}.
\end{align*}

Therefore, $f\leq g$ on $X.$

b) The proof follows immediately from a) because $f(x)=g(x)$ means that both
relations $f(x)\leq g(x)$ and $f(x)\geq g(x)$ are satisfied.
\end{proof}

\section{Dedekind order completeness of $\mathbb{H}_{f\!t}(X)$}

In this section we will prove that $\mathbb{H}_{f\!t}(X)$ is Dedekind order
complete, see Theorems \ref{tocomp1} and \ref{tocomp2}. Here we recall that
the partial order on $\mathbb{H}_{f\!t}(X)$ is the one induced by the
partial order on $\mathbb{A}(X)$. And the partial order on $\mathbb{A}(X)$
is induced by (\ref{iorder}) point-wise on the respective functions similar
with (\ref{order}).

Upon an obvious extension of the respective result in \cite{Baire} we have
the following lemma about semi-continuous functions.

\begin{lemma}
\label{lsupl}a) Let $L\subseteq\mathcal{A}(X)$ be a set of lower
semi-continuous functions. Then the function $l$ defined by $l(x)=\sup
\{f(x):f\in L\}$ is lower semi-continuous.\newline
b) Let $U\subseteq \mathcal{A}(X)$ be a set of upper semi-continuous
functions. Then function the $u$ defined by $u(x)=\inf\{f(x):f\in U\}$ is
upper semi-continuous.
\end{lemma}

\begin{theorem}
\label{tocomp1}Let $\mathcal{F}$ be a subset of $\mathbb{H}_{f\!t}(X)$ which
is bounded from above. Then there exists $u\in\mathbb{H}_{f\!t}(X)$ such
that $u=\sup\mathcal{F}$.
\end{theorem}

\begin{proof}
Let $\psi\in\mathbb{H}_{f\!t}(X)$ be an upper bound of $\mathcal{F}$, that
is, $f\leq\psi$ for all $f\in\mathcal{F}$. Denote
\begin{equation*}
g(x)=\sup\{\underline{f}(x):f=[\underline{f},\overline{f}]\in\mathcal{F}\}%
\text{, }x\in X.
\end{equation*}
It is easy to see that%
\begin{equation}
-\infty<g(x)\leq\psi(x)\text{, }x\in X.  \label{ineq20}
\end{equation}

The function $g$ being a supremum of lower semi-continuous functions is also
a lower semi-continuous function, see Lemma \ref{lsupl}. Therefore,
according to Theorem \ref{FSI} we have $u=F(S(g))=F(S(I(g)))\in\mathbb{H}(X)$%
. Furthermore, from the monotonicity of the operators $S$ and $F$ and the
inequality (\ref{ineq20}) we have%
\begin{align*}
u(x) & =F(S(g))(x)\leq F(S(\psi))(x)=\psi(x)<\infty, \\
u(x) & =F(S(g))(x)\geq I(S(g))(x)\geq I(g)(x)=g(x)>-\infty.
\end{align*}

Therefore, $u\in\mathbb{H}_{f\!t}(X)$. We will prove that $u$ is the
supremum of $\mathcal{F}$. More precisely, we will show that

A) $u$ is an upper bound of $\mathcal{F}$, i.e. $f(x)\leq u(x),\;\;x\in
X,\;\;f\in\mathcal{F}$

B) $u$ is the smallest upper bound of $\mathcal{F}$, i.e. for any function $%
h\in\mathbb{H}_{f\!t}(X)$,%
\begin{equation*}
f(x)\leq h(x),\;\;x\in X,\;\;f\in\mathcal{F}\;\;\Longrightarrow\;\;u(x)\leq
h(x),\;\;x\in X.
\end{equation*}
Using the monotonicity of the operators $S$ and $F,$ for every $f=[%
\underline {f},\overline{f}]\in\mathcal{F}$ we have
\begin{align*}
f(x) & \leq\overline{f}(x)=S(\underline{f})(x)\leq S(g)(x)\text{, }x\in X, \\
f(x) & =F(f)(x)\leq F(S(g))(x)=u(x)\text{, }x\in X,
\end{align*}%
i.e., $u$ is an upper bound of the set $\mathcal{F}$.

Assume that function $h\in\mathbb{H}_{f\!t}(X)$ is such that
\begin{equation*}
f(x)\leq h(x),\;\;x\in X,\;\;f\in\mathcal{F}.
\end{equation*}
We have
\begin{equation*}
g(x)=\sup\{I(f)(x):f\in\mathcal{F} \}\leq h(x)\text{, \ } x\in X,
\end{equation*}
which implies $S(g)\leq S(h).$ Hence $u=F(S(g))\leq F(S(h))=h.$ Therefore, $%
u=\sup\mathcal{F}$.
\end{proof}

In a similar manner we have

\begin{theorem}
\label{tocomp2}Let $\mathcal{F}$ be a subset of $\mathbb{H}_{f\!t}(X)$ which
is bounded from below. Then there exists $v\in\mathbb{H}(X)$ such that $%
v=\inf\mathcal{F}$.$\hfill\square$
\end{theorem}

\section{ Dedekind order completion of $C_{b}(X)$}

From the Dedekind order completeness of $\mathbb{H}_{f\!t}(X)$ it follows
that it contains the Dedekind order completion of its every subset, in
particular the set $C_{b}(X)$. We will show that if $X$ is completely
regular the Dedekind order completion of $C_{b}(X)$ is
\begin{equation*}
\mathbb{H}_{b}(X)=\{f\in\mathbb{H}(X):\exists M\in\mathbb{R}:-M\leq f(x)\leq
M,\text{ }x\in X\}.
\end{equation*}
Since the inclusion
\begin{equation*}
C_{b}(X)\subseteq\mathbb{H}_{b}(X)\subseteq\mathbb{H}_{f\!t}(X)
\end{equation*}
is obvious the statement will be proved through the next two theorems. More
precisely, the Theorem \ref{tHb1} shows that the set $\mathbb{H}_{b}(X)$ is
Dedekind order complete while the Theorem \ref{tHb2} shows that $\mathbb{H}%
_{b}(X)$ is the minimal Dedekind order complete set containing $C_{b}(X)$.

\begin{theorem}
\label{tHb1}Let $\mathcal{F}\subseteq\mathbb{H}_{b}(X)$.\newline
a) If $\mathcal{F}$ is bounded from above then there exists $u\in\mathbb{H}%
_{b}(X)$ such that $u=\sup\mathcal{F}.$\newline
b) If $\mathcal{F}$ is bounded from below then there exists $v\in\mathbb{H}%
_{b}(X)$ such that $v=\inf\mathcal{F}$.
\end{theorem}

\begin{proof}
a) Let $\mathcal{F}$ be bounded from above, that is, there exists a function
$\psi\in\mathbb{H}_{b}(X)$ such that $f\leq\psi$ for all $f\in\mathcal{F}$.
Since $\mathbb{H} _{f\!t}(X)$ is Dedekind order complete the set $\mathcal{F}
$ has a supremum in $\mathbb{H}_{f\!t}(X)$. Let $u=\sup\mathcal{F} \in
\mathbb{H}_{f\!t}(X)$. We need to show that $u\in\mathbb{H}_{b}(X).$ Since $%
\psi$ is an upper bound of $\mathcal{F}$, we have $u=\sup\mathcal{F}\leq\psi
$. Let $f_{0}\in\mathcal{F}$. Then
\begin{equation*}
f_{0}\leq u\leq\psi.
\end{equation*}
Using that both $f_{0}$ and $\psi$ are in $\mathbb{H}_{b}(X)$ there exist
real numbers $M_{1}$ and $M_{2}$ such that
\begin{equation*}
-M_{1}\leq f_{0}(x)\leq M_{1}\text{ , \ }-M_{2}\leq\psi(x)\leq M_{2}\text{ ,
\ }x\in X.
\end{equation*}
Hence, for $M=\max\{M_{1},M_{2}\}$ we have%
\begin{equation*}
-M\leq-M_{1}\leq f_{0}(x)\leq u(x)\leq\psi(x)\leq M_{2}\leq M\text{ , \ }%
x\in X.
\end{equation*}
Therefore, $u\in\mathbb{H}_{b}(X).$

The statement in b) is proved in a similar way.
\end{proof}

\begin{theorem}
\label{tHb2}If the topological space $X$ is completely regular, for every
function $f\in\mathbb{H}_{b}(X)$ there exists a set $\mathcal{F}\subseteq
C(X)$ such that $f=\sup\mathcal{F}$.
\end{theorem}

\begin{proof}
Let $f=[\underline{f},\overline{f}]\in\mathbb{H}_{b}(X)$ be bounded by $M$,
that is,
\begin{equation*}
-M\leq\underline{f}(x)\leq\overline{f}(x)\leq M\text{, \ }x\in X,
\end{equation*}
and let
\begin{equation}
\mathcal{F}=\{g\in C(X):g\leq f\}.  \label{calF}
\end{equation}

We will show that $f=\sup\mathcal{F}$. Let $\varepsilon$ be an arbitrary
positive number and let $x\in X$. Since the function $\underline{f}$ is
lower semi-continuous, there exists a neighborhood $V_{x}\in\mathcal{V}_{x}$
such that
\begin{equation}
\underline{f}(y)>\underline{f}(x)-\varepsilon,\text{ }y\in V_{x}.
\label{ineq12}
\end{equation}

Denote $m_{x}=\inf_{y\in V_{x}}\underline{f}(y)$. We have%
\begin{equation}
\underline{f}(x)-\varepsilon\leq m_{x}\leq\underline{f}(y),\text{ } y\in
V_{x}.  \label{ineq13}
\end{equation}

Due to the complete regularity of $X$ there exists a function $\phi_{x}$
such that
\begin{align}
\phi_{x}(x) & =1,  \notag \\
\phi_{x}(y) & =0\text{ \ if \ }y\notin V_{x},  \label{propreg} \\
0 & \leq\phi_{x}(y)\leq1\text{ for all } y\in X.  \notag
\end{align}

Consider the function $g_{x}$ defined by%
\begin{equation*}
g_{x}(y)=(m_{x}+M)\phi_{x}(y)-M,\text{ }y\in X.
\end{equation*}
Since $\phi_{x}$ is continuous then $g_{x}$ is also continuous. Furthermore,
using inequalities (\ref{ineq13}) and the properties (\ref{propreg}) we have%
\begin{align*}
g_{x}(y) & =(m_{x}+M)\times0-M=-M\leq f(y)\text{ \ if \ }y\in X\backslash
V_{x}, \\
g_{x}(y) & \leq(m_{x}+M)\times1-M=m_{x}\leq f(y)\text{ \ if \ }y\in V_{x}.
\end{align*}

Therefore, for every $x\in X$ the function $g_{x}$ belongs to the set $%
\mathcal{F}$. Hence%
\begin{equation*}
(\sup\mathcal{F})(x)\geq g_{x}(x)=(m_{x}+M)\phi_{x}(x)-M=m_{x}\geq
\underline{f}(x)-\varepsilon,\text{ }x\in X,
\end{equation*}
where the last inequality follows from (\ref{ineq13}). Since the positive $%
\varepsilon$ in the above inequality is arbitrary we have%
\begin{equation*}
(\sup\mathcal{F})(x)\geq\underline{f}(x),\text{ }x\in X.
\end{equation*}
Using the monotonicity of the operator $F$ and the fact that both $\sup%
\mathcal{F}$ and $f$ are H-continuous we obtain%
\begin{equation}
\sup\mathcal{F}\geq F(\underline{f})=f.  \label{ineq15}
\end{equation}

However, $f$ is an upper bound of $\mathcal{F}$, see (\ref{calF}), while $%
\sup\mathcal{F}$ is the smallest upper bound of $\mathcal{F}$ in $\mathbb{H}%
_{b}(X)$. Thus, from the inequality (\ref{ineq15}) it follows that $\sup%
\mathcal{F}=f$ which completes the proof.
\end{proof}

The Dedekind order completion $C_{b}(X)^{\#}=\mathbb{H}_{b}(X)$ of $C_{b}(X)$
discussed in this section is similar to the result of Dilworth in \cite%
{Dilworth} where it is proved that $C_{b}(X)^{\#}$ is order isomorphic to
the set of the so called normal upper semi-continuous functions. The
approach presented here produces an alternative characterization of $%
C_{b}(X)^{\#}$ as a set of Hausdorff continuous functions. The power if this
approach is further demonstrated in the next section by constructing the
Dedekind order completion of the set $C(X)$ without any boundedness
assumption about the functions, a case which is particularly important for
applications to generalized solutions of nonlinear partial differential
equations.

\section{Dedekind order completion of $C(X)$}

Since $C(X)$ is a subset of $\mathbb{H}_{f\!t}(X)$, see (\ref{incl1}), which
is Dedekind order complete, see Theorems \ref{tocomp1} and \ref{tocomp2},
the Dedekind order completion of $C(X)$ is a subset of $\mathbb{H}_{f\!t}(X)$
for any topological space $X$. We will show that for a completely regular
topological space $X$ the Dedekind order completion of $C(X)$ is the set $%
\mathbb{H}_{cm}(X)$ of all Hausdorff continuous functions with a continuous
majorant and a continuous minorant. More precisely,%
\begin{equation*}
\mathbb{H}_{cm}(X)=\{f\in\mathbb{H}_{f\!t}(X):\exists\varphi,\psi\in
C(X):\varphi\leq f\leq\psi\}
\end{equation*}
and we will prove that
\begin{equation}
C(X)^{\#}=\mathbb{H}_{cm}(X).  \label{cxHcm}
\end{equation}
We have the inclusion $C(X)\subseteq\mathbb{H}_{cm}(X)$. Similarly to the
preceding section we will prove (\ref{cxHcm}) by showing that $\mathbb{H}%
_{cm}(X)$ is Dedekind order complete, see Theorem \ref{tHcm1}, and that $%
\mathbb{H}_{cm}(X)$ is the minimal Dedekind order complete set which
contains $C(X)$, see Theorem \ref{tHcm2}.

\begin{theorem}
\label{tHcm1}Let $\mathcal{F}\subseteq\mathbb{H}_{cm}(X)$.\newline
a) If $\mathcal{F}$ is bounded from above then there exists $u\in\mathbb{H}%
_{cm}(X)$ such that $u=\sup\mathcal{F}$.\newline
b) If $\mathcal{F}$ is bounded from below then there exists $v\in\mathbb{H}%
_{cm}(X)$ such that $v=\inf\mathcal{F}$.
\end{theorem}

\begin{proof}
a) Let $\mathcal{F}$ be bounded from above, that is, there exists a function
$h\in\mathbb{H}_{cm}(X)$ such that $f\leq h$ for all $f\in\mathcal{F}$.
Since $\mathbb{H}_{f\!t}(X)$ is Dedekind order complete the set $\mathcal{F}$
has a supremum in $\mathbb{H}_{f\!t}(X)$. Denote $u=\sup\mathcal{F}\in
\mathbb{H}_{f\!t}(X)$. We need to show that $u\in\mathbb{H}_{cm}(X).$ Since $%
h$ is an upper bound of $\mathcal{F}$, we have $u=\sup\mathcal{F}\leq h$.
Let $f_{0}\in\mathcal{F}$. Then%
\begin{equation*}
f_{0}\leq u\leq h.
\end{equation*}
Using that both $f_{0}$ and $h$ are in $\mathbb{H}_{cm}(X)$ there exist real
functions $\varphi_{1},\varphi_{2},\psi_{1},\psi_{2}\in C(X)$ such that
\begin{equation*}
\varphi_{1}\leq f_{0}\leq\psi_{1}\text{ , \ }\varphi_{2}\leq h\leq\psi_{2}.
\end{equation*}
Then, we have%
\begin{equation*}
\varphi_{1}\leq f_{0}\leq u\leq h\leq\psi_{2}.
\end{equation*}
Hence, $u\in\mathbb{H}_{cm}(X).$

The statement in b) is proved in a similar way.
\end{proof}

\begin{theorem}
\label{tHcm2}If $X$ is a completely regular topological space for every
function $f\in\mathbb{H}_{cm}(X)$ there exists a set $\mathcal{F}\subseteq
C(X)$ such that $f=\sup\mathcal{F}$.
\end{theorem}

\begin{proof}
Let $f=[\underline{f},\overline{f}]\in\mathbb{H}_{cm}(X)$ and let $%
\varphi\in C(X)$ be a minorant of $f$, that is $\varphi\leq f$. Consider the
set
\begin{equation}
\mathcal{F}=\{g\in C(X):g\leq f\}.  \label{calF2}
\end{equation}
Clearly, $\varphi\in\mathcal{F}$. Hence the set $\mathcal{F}$ is not empty.\
We will show that $f=\sup\mathcal{F}$. Let $\varepsilon$ be an arbitrary
positive number and let $x\in X$. Since the function $\underline{f}$ is
lower semi-continuous, there exists a neighborhood $V_{x} \in\mathcal{V}_{x}$
such that
\begin{equation}
\underline{f}(y)>\underline{f}(x)-\varepsilon,\text{ }y\in V_{x}.
\label{ineq122}
\end{equation}
Denote $m_{x}=\inf_{y\in V_{x}} \underline{f}(y)$. From (\ref{ineq122}) it
follows that
\begin{equation}
\underline{f}(x)-\varepsilon\leq m_{x}\leq\underline{f}(y),\text{ }y\in
V_{x}.  \label{ineq132}
\end{equation}
Since the space $X$ is completely regular there exists a function $%
\phi_{x}\in C(X)$ such that
\begin{align}
\phi_{x}(x) & =1,  \notag \\
\phi_{x}(y) & =0\text{ \ if \ }y\notin V_{x},  \label{propreg2} \\
0 & \leq\phi_{x}(y)\leq1\text{ for all }y\in X.  \notag
\end{align}
Consider the function $g_{x}$ defined by
\begin{equation*}
g_{x}(y)=(m_{x}-\varphi(y))\phi_{x}(y)+\varphi(y),\text{ }y\in X.
\end{equation*}
Since $\phi_{x},\varphi\in C(X)$ then $g_{x}\in C(X)$. Furthermore, using
inequalities (\ref{ineq132}) and the properties (\ref{propreg2}) we have
\begin{align*}
g_{x}(y) & =(m_{x} -\varphi(y))\times0+\varphi(y)=\varphi(y)\leq f(y)\text{
\ if \ } y\in X\backslash V_{x}\text{ ,} \\
g_{x}(y) & \leq0\times\phi_{x} (y)+\varphi(y)=\varphi(y)\leq f(y)\text{ \ if
}\varphi(y)>m_{x} \text{ ,} \\
g_{x}(y) & \leq(m_{x}-\varphi(y))\times1+\varphi(y)=m_{x} \leq f(y)\text{ \
if \ } y\in V_{x}\text{ and }\varphi(y)<m_{x} \text{ .}
\end{align*}
Therefore, for every $x\in X$ we have $g_{x}\leq f$ which implies $g_{x}\in
\mathcal{F}$. Hence
\begin{equation*}
(\sup\mathcal{F})(x)\geq g_{x}(x)=(m_{x}-\varphi(x))\varphi_{x}(x)+\varphi
(x)=m_{x}\geq\underline{f} (x)-\varepsilon,\text{ } x\in X,
\end{equation*}
where the last inequality follows from (\ref{ineq132}). Since the positive $%
\varepsilon$ in the above inequality is arbitrary we have
\begin{equation*}
(\sup\mathcal{F})(x)\geq\underline{f}(x),\text{ } x\in X.
\end{equation*}
Using the monotonicity of the operator $F$ and the fact that both $\sup%
\mathcal{F}$ and $f$ are H-continuous we obtain
\begin{equation}
\sup\mathcal{F}\geq F(\underline{f})=f.  \label{ineq152}
\end{equation}
However, $f$ is an upper bound of $\mathcal{F}$, see (\ref{calF2}), while $%
\sup\mathcal{F}$ is the smallest upper bound of $\mathcal{F}$ in $\mathbb{H}%
_{cm}(X)$. Thus, from the inequality (\ref{ineq152}) it follows that $\sup%
\mathcal{F}=f$ which completes the proof.
\end{proof}

The characterizations of the Dedekind order completion of $C(X)$ as the set $%
\mathbb{H}_{cm}(X)$ is interesting if the topological space $X$ is such that
$C(X)$ contains unbounded functions because the case of bounded continuous
functions was dealt with in section 4. For example, if $X$ is a compact
Hausdorff space, the elements of both $C(X)$ and $\mathbb{H}_{f\!t}(X)$ are
bounded functions, that is, we have%
\begin{equation*}
C(X)=C_{b}(X)\text{ and }\mathbb{H}_{b}(X)=\mathbb{H}_{f\!t}(X),
\end{equation*}
which means that the Dedekind order completion of $C(X)$ is $\mathbb{H}%
_{f\!t}(X)$. Since $\mathbb{H}_{cm}(X)=C(X)^{\#}$, this also means that
\begin{equation}
\mathbb{H}_{cm}(X)=\mathbb{H}_{f\!t}(X).  \label{HcmHft}
\end{equation}
The characterization of the topological spaces for which the condition (\ref%
{HcmHft}) holds is still an open problem. However, we can show that (\ref%
{HcmHft})\ holds if $X$ is an arbitrary metric space. Note that in this case
$C(X)$ may contain unbounded functions.

\begin{theorem}
Let $(X,\rho)$ be a metric space. Then, every function $f\in\mathbb{H}%
_{f\!t}(X)$ has a continuous majorant and a continuous minorant, that is, $%
\mathbb{H}_{cm}(X)=\mathbb{H}_{f\!t}(X)$.
\end{theorem}

\begin{proof}
We will use the function $h:\mathbb{R}\rightarrow (-1,1)\subset R$ defined
by
\begin{equation*}
h(z)=\frac{z}{1+|z|},\text{ }z\in \mathbb{R}\text{.}
\end{equation*}%
This real function is continuous and strictly increasing. The inverse
function $h^{-1}:(-1,1)\rightarrow \mathbb{R}$ is given by
\begin{equation*}
h^{-1}(z)=\frac{z}{1-|z|},\text{ }z\in (-1,1),
\end{equation*}%
and is also continuous and strictly increasing.

Let $f=[\underline{f},\overline{f}]\in\mathbb{H}_{f\!t}(X)$. Consider the
function $\varphi:X\times X\rightarrow R$ defined by%
\begin{equation}
\varphi(t,x)=h(\underline{f}(t))+\rho(t,x)=\frac{\underline{f}(t)} {1+|%
\underline{f}(t)|}+\rho(t,x).  \label{fhi}
\end{equation}
It is easy to see that the function $\varphi$ is bounded from below. Indeed,
since the value of the metric $\rho$ is always nonnegative and the fraction
in (\ref{fhi}) is greater than -1 we have $\varphi(t,x)>-1$. Then, we can
define
\begin{equation*}
\psi(x)=\inf\{\varphi(t,x):t\in X\}.
\end{equation*}
First we will show that the function $\psi$ is continuous on $X$. From the
triangular inequality of the metric $\rho$ for every $x,y,t\in X$ we have%
\begin{equation*}
\rho(t,y)-\rho(x,y)\leq\rho(t,x)\leq\rho(t,y)+\rho(x,y).
\end{equation*}
Therefore%
\begin{equation*}
\varphi(t,y)-\rho(x,y)\leq\varphi(t,x)\leq\varphi(t,y)+\rho(x,y).
\end{equation*}
Taking the infimum on $t\in X$ we obtain
\begin{equation*}
\psi(y)-\rho(x,y)\leq\psi(x)\leq\psi(y)+\rho(x,y).
\end{equation*}
Hence we have the inequality%
\begin{equation*}
|\psi(x)-\psi(y)|\leq\rho(x,y),
\end{equation*}
which implies that the function $\psi$ is continuous on $X$.

Our second step is to prove that $\psi$\ satisfies the inequalities%
\begin{equation}
-1<\psi(x)\leq\frac{\underline{f}(x)}{1+|\underline{f}(x)|}<1,\text{ } x\in
X.  \label{ineq19}
\end{equation}
For every $x\in X$ we have
\begin{equation*}
\psi(x)=\inf\{\varphi(t,x):t\in X\}\leq\varphi(x,x)=\frac{\underline{f}(x)}{%
1+|\underline{f}(x)|}.
\end{equation*}
Furthermore, since $-1$ is a lower bound of $\varphi(t,x)$ the inequality
\begin{equation*}
\psi(x)\geq-1
\end{equation*}
also holds. It remains to prove that $\psi(x)\neq-1$. Let us assume that
there exists $x\in X$ such that $\psi(x)=-1$. Let the real number $\mu$ be
such that $-1<\mu<\frac{\underline{f}(x)}{1+|\underline{f}(x)|}$. Then we
have%
\begin{equation*}
h(\underline{f}(x))=\frac{\underline{f}(x)}{1+|\underline{f}(x)|} >\mu>-1.
\end{equation*}
Due to the continuity of $h$ there exists $\eta$ such that
\begin{equation*}
h(z)>\mu\text{ whenever } \underline{f}(x)-\eta<z<\underline{f}(x)+\eta.
\end{equation*}
Moreover, since $h$ is strictly increasing we have
\begin{equation*}
h(z)>\mu\text{ for all }z>\underline{f}(x)-\eta.
\end{equation*}
Using that the function $\underline{f}(t)$ is lower semi-continuous at $x$,
there exists $\varepsilon>0$ such that
\begin{equation*}
\underline{f}(t)>\underline{f}(x)-\eta\text{ \ whenever \ }%
\rho(t,x)<\varepsilon.
\end{equation*}
Hence
\begin{equation}
\frac{\underline{f}(t)}{1+|\underline{f}(t)|}=h(\underline{f}(t))>\mu\text{
\ whenever \ } \rho(t,x)<\varepsilon.  \label{ineq16}
\end{equation}

Let now $\delta=\min\{\varepsilon,\mu+1\}$. Since $\psi(x)$ is defined as an
infimum on $t\in X$, there exists $t_{\delta} \in X$ such that
\begin{equation*}
-1=\psi(x)\leq\varphi(t_{\delta},x)\leq\psi(x)+\delta=-1+\delta
\end{equation*}
or, more precisely,
\begin{equation*}
-1\leq\frac{\underline{f} (t_{\delta})}{1+|\underline{f}(t_{\delta})|}%
+\rho(t_{\delta} ,x)\leq-1+\delta.
\end{equation*}
Using simple manipulations we obtain%
\begin{align}
0 & \leq\rho(t_{\delta},x)\leq\delta-\left( 1+\frac{\underline{f}(t_{\delta
})}{1+|\underline{f}(t_{\delta})|} \right) <\delta\leq\varepsilon
\label{ineq17} \\
-1 & \leq\frac{\underline{f}(t_{\delta})}{1+|\underline{f} (t_{\delta})|}%
\leq-1+\delta\leq\mu  \label{ineq18}
\end{align}

The contradiction between inequalities (\ref{ineq17}), (\ref{ineq18}) on the
one side and the condition (\ref{ineq16}) on the other side show that the
assumption that $\psi(x)=-1$ for some $x\in X$ is false. Therefore $%
\psi(x)>-1,$ $x\in X$.

Now we can consider the function%
\begin{equation*}
\phi(x)=h^{-1}(\psi(x))=\frac{\psi(x)}{1-|\psi(x)|},\text{ }x\in X.
\end{equation*}
Due to inequalities (\ref{ineq19}) function $\phi$ is well defined for every
$x\in X$. Furthermore, $\phi$ is continuous on $X$ because $\psi$ is
continuous on $X$. Using the fact that the function $h^{-1}$ is strictly
increasing on the interval $(-1,1)$ from the middle inequality in (\ref%
{ineq19}) we obtain%
\begin{equation*}
\phi(x)=h^{-1}(\psi(x))\leq h^{-1}\left( \frac{\underline{f}(x)}{1+|%
\underline{f}(x)|}\right) =h^{-1}\left( h(\underline{f}(x)\right) =%
\underline{f}(x)\leq f(x).
\end{equation*}
Thus, $\phi$ is a continuous minorant of $f$. The existence of a continuous
majorant is proved in a similar way.
\end{proof}

\begin{corollary}
If $X$ is a metric space then $C(X)^{\#}=\mathbb{H}_{f\!t}(X).$
\end{corollary}

\section{Conclusion}

This paper addresses the issue of finding a characterization for the
Dedekind order completion of $C(X)$ in a constructive form, namely, as a set
of functions on the same space $X$. Here $X$ can be an arbitrary topological
space. The functions, which give the completion, are the Hausdorff
continuous functions, and as such, are in general interval valued functions.
However, they are not unlike the usual real-valued functions, since interval
values only appear at points of discontinuity. The set $\mathbb{H}_{f\!t}(X)$
of all finite Hausdorff continuous functions is Dedekind order complete and,
thus, it contains the Dedekind completion of $C(X)$. It was shown that in
the important case of $X$ being a metric space $\mathbb{H}_{f\!t}(X)$ is the
Dedekind order completion of $C(X)$. In the more general case of a
completely regular topological space $X$ the Dedekind order completion of $%
C(X)$ is the set $\mathbb{H}_{cm}(X)$ of all H-continuous functions with a
continuous majorant and a continuous minorant. In the case considered by
Dilworth \cite{Dilworth}, namely, the Dedekind order completion of $C_{b}(X)$%
, where $X$ is a completely regular space, an alternative characterization
of the Dedekind order completion of $C_{b}(X)$ was found in the form of the
set $\mathbb{H}_{b}(X)$ - the set of all bounded Hausdorff continuous
functions.

The Hausdorff continuous functions can be applied in further improving
recent mathematical methods in the solution of large classes on nonlinear
partial differential equations, methods which use the Dedekind order
completion of $C(X)$, \cite{Rosinger}. Let us note that the applicability of
the set $\mathbb{H}_{f\!t}(X)$ in the order completion method discussed in %
\cite{Rosinger} depends on the fact that, unlike the result in \cite%
{Dilworth}, there is no boundedness requirement on the functions. In
addition to that, the Hausdorff continuous functions provide a convenient
representation of nonlinear shock waves. For instance, the function in
example \ref{shock}, as mentioned, is a solution of the inviscid Burger's,
or nonlinear shock wave equation, with a shock along the line $t=1+2x,x\geq
0 $.

\section{Acknowledgements}

The author would like to thank Prof. E E Rosinger for his encouragement,
continued support and useful discussions.

\end{document}